\documentclass[journal]{IEEEtran}
\usepackage{amsfonts}
\usepackage{amssymb,amsmath,color,graphicx}
\usepackage{subfigure}
\usepackage{cite}
\usepackage{comment}
\usepackage{verbatim}
\usepackage{algorithm}
\usepackage{algorithmic}
\usepackage{xcolor}
\usepackage{placeins}
\usepackage{color,soul}
\usepackage{setspace}
\usepackage{float}
\usepackage{bm}
\usepackage{dsfont}
\usepackage{xr}
\usepackage{tabu}
\usepackage{tabularx}
\usepackage{longtable}
\usepackage[utf8]{inputenc}
\usepackage{pifont}
\usepackage{wasysym}
\usepackage{lipsum}
\usepackage{booktabs}
\usepackage{tabulary}
\usepackage{textcomp}
\usepackage{booktabs}
\usepackage{multirow}
\usepackage{multicol}
\usepackage{wrapfig}
\usepackage{mathtools}
\usepackage{subfig}
\usepackage{booktabs}
\usepackage{array, makecell}
\usepackage{caption}
\usepackage{subcaption}
\allowdisplaybreaks

\ifCLASSINFOpdf
\else
\fi

\begin{document}
% Titles are generally capitalized except for words such as a, an, and, as,
% at, but, by, for, in, nor, of, on, or, the, to and up, which are usually
% not capitalized unless they are the first or last word of the title.
\title{Spatial Arbitrage Through Bidirectional Electric Vehicle Charging with Delivery Fleets}

\author{Mostafa~Mohammadian,~\IEEEmembership{Student~Member,~IEEE,}
        Constance~Crozier,~\IEEEmembership{Member,~IEEE,}\\
        ~Kyri~Baker,~\IEEEmembership{Senior Member,~IEEE}% <-this % stops a space
\thanks{M. Mohammadian and K. Baker are with the University of Colorado Boulder, Boulder, CO $80309$ USA. (email:mostafa.mohammadian@colorado.edu; kyri.baker@colorado.edu). C. Crozier is with Georgia Institute of Technology, Atlanta, GA $30332$ (email: ccrozier8@gatech.edu).}% <-this % stops a space
% \thanks{This work was supported by the National Science Foundation CAREER award 2041835. Additionally, this work utilized the Summit supercomputer, which is supported by the National Science Foundation (awards ACI-$1532235$ and ACI-$1532236$), the University of Colorado Boulder, and Colorado State University. The Summit supercomputer is a joint effort of the University of Colorado Boulder and Colorado State University.}
}

\maketitle

\begin{abstract}
The adoption of electric vehicles (EVs), including electric taxis and buses, as a mode of transportation, is rapidly increasing in cities. In addition to providing economic and environmental benefits, these fleets can potentially participate in the energy arbitrage market by leveraging their mobile energy storage capabilities. This presents an opportunity for EV owners to contribute to a more sustainable and efficient energy system while also reducing their operational costs. The present study introduces deterministic and single-stage stochastic optimization frameworks that aim to maximize revenue by optimizing the charging, discharging, and travel of a fleet of electric vehicles in the context of uncertainty surrounding both spatial and temporal energy prices.
The simulations are performed on a fleet of electric delivery trucks, which have to make deliveries to certain locations on specific dates.
The findings indicate the promising potential of bidirectional electric vehicle charging as a mobile grid asset. However, it is important to note that significant revenue is only realized in scenarios where there is substantial variation in prices between different areas, and when these price variations can be accurately forecasted with a high level of confidence.
\end{abstract}


\IEEEpeerreviewmaketitle

\section*{Nomenclature}
\addcontentsline{toc}{section}{Nomenclature}
\begin{IEEEdescription}[\IEEEusemathlabelsep\IEEEsetlabelwidth{$V_1,V_2,V_3$}]
\item [\textbf{Parameters}]
\item[$N$] Total number of vehicles
\item[$T$]  Total number of time steps
\item[$t$] Discrete time slot index
\item[$\Delta t$] Size of a timestep (hours)
\item[$E_{n, 0}$] Initial energy in the battery of $n^\text{th}$ vehicle 
\item[$E^{\max}_n$] Maximum battery energy of $n^\text{th}$ vehicle
\item[$P^{drive}_n$] Discharge rate of $n^\text{th}$ vehicle when driving
\item[$P_{d}^{\max}$] Maximum discharge rate
\item[$P_{c}^{\max}$] Maximum charge rate
\item[$T_{A,B}$] Number of time steps takes to travel from A to B (or reverse)
\item[$T_{A,C}$] Number of time steps takes to travel from A to C (or reverse)
\item[$T_{B,C}$] Number of time steps takes to travel from B to C (or reverse)
\item[$\eta_{c}$] Charging efficiency
\item[$\eta_{d}$] Discharging efficiency

\item [\textbf{Variables}]
\item[$x,y,z$] EV’s charged (positive) or discharged (negative) power at locations
$A$, $B$, and $C$
\item[$\makecell{\alpha_{n, t}, \beta_{n, t},\\ \newline \gamma_{n, t}}$] (Binary) symbolizing the presence of $n^\text{th}$ vehicle is at locations $A$, $B$, and $C$
% \item $\beta_t^{(j)}$ (binary) symbolizes that the vehicle $j$ is at location B
% \item $\gamma_t^{(j)}$ (binary) symbolizes that the vehicle $j$ is at location C
\item[$c_{n, t}^A$] Grid-side charging power of $n^\text{th}$ vehicle at location $A$
\item[$d_{n, t}^A$] Grid-side discharging power of $n^\text{th}$ vehicle at location $A$
\item[$c_{n, t}^B$] Grid-side charging power of $n^\text{th}$ vehicle at location $B$
\item[$d_{n, t}^B$] Grid-side discharging power of $n^\text{th}$ vehicle at location $B$
\item[$c_{n, t}^C$] Grid-side charging power of $n^\text{th}$ vehicle at location $C$
\item[$d_{n, t}^C$] Grid-side discharging power of $n^\text{th}$ vehicle at location $C$
\item[$\makecell{a_{n}, b_n,\\ c_n}$] (Binary) symbolizing the $n^\text{th}$ vehicle being at location A, B, and C at any time instance.
\end{IEEEdescription}

\section{Introduction}
\IEEEPARstart{W}{ith} the advent of distributed energy resource capabilities and new market mechanisms, the energy arbitrage (e.g., buying during “off-peak” times and selling during “on-peak” times) landscape is rapidly changing from large-scale stationary energy storage owned by large firms to mobile energy storage owned by individuals with electric vehicles (EVs) \cite{7892020}. EVs have the potential to impact energy arbitrage in two significant ways: first, the mobility of EVs allows for them to act as distributed energy resources (DERs), providing additional flexibility to the grid \cite{Yuan22}; second, the rise of vehicle-to-grid (V2G) technology allows EVs to not only consume energy but also supply excess energy back to the grid \cite{HE2021379}. This shift towards mobile energy storage provides new opportunities for individual EV owners to participate in energy arbitrage and contribute to a more sustainable energy future, due to several major advantages over
traditional gas-powered vehicles such as being environmentally friendly, cheaper to
maintain, and generally safer \cite{10.1145/3538637.3538870}.

The adoption of vehicle-to-grid (V2G) technologies have been rapidly increasing, facilitated by the proliferation of EVs in the energy system \cite{10102467}. This trend aligns with the projected growth of EVs, with an estimated $30–42$ million EVs expected to be on the roads by $2030$ in the U.S. \cite{NREL}. To support this transition, significant advancements in bidirectional charging have been achieved. Major automotive manufacturers like Ford are introducing electric options for their top-selling vehicles, highlighting the expanding capabilities of bidirectional charging. In a significant regulatory development, the Federal Energy Regulatory Commission (FERC) approved Order $2222$ in $2020$, allowing distributed energy resources, including EVs, to participate in wholesale markets \cite{FERC}. These advancements not only facilitate the integration of EVs into the energy system but also open up new opportunities for their contribution to the overall energy landscape.

The scheduling of EVs for profit maximization is challenging due to uncertainties in EV mobility behavior and market price volatility. This complexity arises in the decision-making process under uncertainty when optimizing EV scheduling \cite{9737306,Sarker16}. Previous studies have demonstrated the cost-effectiveness of V2G technologies in certain situations, while the feasibility of the V2G-enabled fleet business model depends on evaluating the financial risks associated with accelerated battery degradation \cite{en13184742}. 
Moreover, the advantages of V2G extend to various scales, including the transmission network, distribution network, and microgrids. By adopting an optimal bidirectional V2G mode for EVs, both the operational stability and economic performance of the system can be enhanced \cite{7802551, CROZIER2020114214}.

Several studies have put forth optimization strategies aimed at maximizing the profit of EVs across different market scenarios. These strategies encompass a range of applications, such as providing frequency regulation \cite{10174660}, operating costs and $\text{CO}_2$ emissions \cite{AHMADI2023109195,8973708}, trading in day-ahead and intra-day electricity markets \cite{ZHANG2023121063}, local energy trading \cite{NIZAMI2020114322}, and minimizing grid fluctuations \cite{10175218}. Additionally, some studies have explored the potential for EVs to participate in multiple grid services \cite{10134774}. However, it is important to note that these approaches typically assume a fixed location for the EV, disregarding the fact that electricity prices can exhibit significant variations over short distances. Considering this aspect could unlock further potential benefits in optimizing EV operations. Therefore, unlike other V2G optimization studies that use fixed EV locations \cite{8951120,10274855,7463483}, our model uniquely captures the value of mobility in spatial arbitrage.
%For instance, in Fig. \ref{fig:ERCOT}, a snapshot of the nodal locational marginal price (LMP) map from the Electric Reliability Council of Texas (ERCOT) is presented, highlighting the significant spatial variation in LMP values. The map showcases negative LMPs in close proximity to LMPs reaching hundreds of dollars, indicating the presence of substantial spatial disparities. This observation suggests the potential for capturing additional owner benefits through the integration of spatial characteristics.
% \begin{figure}[h!]
% \centering
% % \includegraphics[clip,trim={0cm 0cm 0cm 1cm},width=\columnwidth]{figs/IMG_1786.jpg}%
% \includegraphics[clip,trim={0cm 0cm 0cm 0cm},width=\columnwidth]{figs/ERCOT.pdf}\\\vspace{-1mm}%
% \caption{LMPs across ERCOT showing negative/very low prices in southern
% Texas and prices exceeding $500\$$/MWh near Austin. ERCOT can have real
% time prices up to and exceeding $9,000\$$/MWh. Source: \cite{ERCOT123}}
% \label{fig:ERCOT}
% \end{figure}

In areas with locational marginal prices (LMPs), time-varying price signals can help identify areas in need of congestion relief and better represent the cost to deliver power to specific locations at specific times. New opportunities for mobile energy storage arise in these situations, capturing both spatial and temporal benefits, which we aim to explore in this paper.

Previous studies have investigated the use of mobile energy storage at the transmission level, focusing on spatio-temporal arbitrage using utility-scale batteries \cite{HE2021379}. However, these studies were limited by a deterministic problem formulation and a restricted driving radius of 10 miles. Therefore, it is essential to further analyze the potential benefits of spatio-temporal LMP differences in a more comprehensive and realistic manner. To leverage spatial price differences, several studies have introduced optimization methods for vehicle charging that incorporate joint optimization of vehicle routing \cite{9172098, 8126871}. These approaches take into account factors such as traffic conditions and charging location uncertainties. However, they do not consider the potential for arbitrage, limiting their focus to unidirectional charging. Furthermore, the complexity of these methods is significantly high due to the extensive number of possible routes. In order to address these limitations, further enhancements are required to incorporate bidirectional charging and mitigate the computational complexity associated with optimizing vehicle routing.

This work represents an expanded and enhanced version of our preliminary conference paper \cite{9916944}. The present paper aims to explore the potential of spatial (also called geographical) arbitrage in the context of bidirectional EV charging given more practical vehicle constraints such as delivery schedules. By combining the objective of arbitrage with the EV's role as a mobile energy storage device, our study focuses on analyzing the potential for fleets of electric delivery trucks to align delivery objectives with a secondary revenue stream resulting from intelligent charging and discharging.

Here, we present a novel approach for addressing the scheduling and optimization challenges faced by a delivery company operating electric delivery trucks between three pre-defined geographical locations. Our proposed framework integrates charging, discharging, delivery, and departure time scheduling into a single-stage stochastic optimization model. To efficiently solve this problem, we formulate it as a stochastic mixed-integer linear program, which can be effectively solved using existing optimization solvers. This framework offers a comprehensive solution for optimizing electric delivery truck operations in a dynamic and uncertain environment.
We conduct a case study in the ERCOT real-time market to evaluate the performance of our proposed framework for electric delivery trucks traveling between San Marcos, San Antonio, and Austin, Texas. The case study allows us to have an initial understanding of the effectiveness and efficiency of the scheduling and optimization strategies in a market setting. Initial results demonstrate that delivery schedules can be adhered to while grid services are provided that benefit both EV owners and the power grid.

% The remainder of this paper is organized as follows. In Section \ref{sec: methods}, after introducing the physics-informed neural networks architecture, we present the extension to g-PINN. Section \ref{sec: model} briefly describes the employed power system model and gPINN used for the prediction task. In Section \ref{sec: result}, we present the simulation results demonstrating the performance of physics-informed neural networks and gPINN methods. Finally, conclusions and future work are given in Section \ref{sec: conclude}.

\section{Methods} \label{sec: methods}
In this section, a stochastic program is introduced for the energy management of a fleet of EVs or electric delivery trucks traveling between three locations. The program is single-stage and employs mixed integer linear optimization techniques. It also considers the charging and discharging behavior of EVs as well as the departure of every single EV, taking into account price and travel uncertainties. We aim to minimize the total charging cost of EVs owned by an aggregator.

This preliminary study is based on several assumptions. Firstly, there is no requirement for EVs to travel separately from their delivery routes; the only requirement is that the EVs make the required delivery by the required time. Additionally, the option for temporal arbitrage is just one subset of the optimization problem. To simplify matters, a separate charging and discharging variable for each vehicle is assumed (e.g. the car can only be charging \emph{or} discharging at any given time). The details of the charging infrastructure, such as the number, type, and locations of chargers, are not explicitly modeled. Instead, maximum charge and discharge rates are specified. The model assumes perfect foresight of real-time electricity prices across all locations and time periods; however, the impact of forecast errors will be analyzed in subsequent sections. It should be noted that this framework has the potential to be extended to include autonomous vehicles that have additional onboard energy storage. Future research could also take into account more complex battery management, degradation effects, and charging strategies, which are important considerations for the economic and technical impacts of V2G technology \cite{9493733}. 

\subsubsection{Electric Vehicle Fleets}
Consider a set $\mathcal{N} = \{1, \ldots, N\}$ of EVs that are going to be charged in different locations $[A, B, C]$ during a specific period of time. To establish a framework for the operation cycle, we partition the horizon time into a set $\mathcal{T} = \{1, \ldots, T\}$ consisting of $T$ discrete time slots. Each of these time slots corresponds to a charging/discharging control interval, typically lasting $15$ minutes. The control interval between any two consecutive slots is represented by $\Delta t$. The assumption made in this case does not impact the general applicability of the proposed algorithm, and it can be adjusted to suit the specific needs of the study case; for example, in markets with different real-time pricing intervals.

The power demand of each EV $n \in \mathcal{N}$ in time slot $t \in \mathcal{T}$ at locations $A, B, C$ has limits $- P_{d}^{\max} < 0$ and $P_{c}^{\max} > 0$, respectively. These limits can reflect the charging and discharging limits of a particular charging station. We thus have the constraints
\begin{subequations}\label{eq: ch/dis variable}
\begin{align}
        - P_{d}^{\max} \leq x_{n, t} \leq P_{c}^{\max}\\
    - P_{d}^{\max} \leq y_{n, t} \leq P_{c}^{\max}\\
    - P_{d}^{\max} \leq z_{n, t} \leq P_{c}^{\max}
\end{align}
\end{subequations}
Note that variables $x_{n, t}, y_{n, t}, z_{n, t}$ represent the $n^{th}$ EV’s charged (positive) or discharged (negative) power at location $A$, $B$, and $C$ at time slot $t \in \mathcal{T}$, respectively. It is worth noting that the charging rate can be influenced by various factors, including the type of charger, battery size, the EV model, and even the current state of charge of the battery. 

Let $10\% \leq E_n^{\text{init}} \leq 90\%$ and $E^{\max}_n \geq 0$ denote the
initial state of charge (SOC) and the maximum capacity of
the battery of each EV $n \in \mathcal{N}$, respectively. Let  $0 < \eta_n^c \leq 1$ and $0 < \eta_n^d \leq 1$ define as the energy transfer efficiencies for charging and discharging of EVs $n \in \mathcal{N}$, respectively.
We introduce slack variables $c_{n, t}^{\{.\}}$ and $d_{n, t}^{\{.\}}$ to indicate the
charging and discharging power of each EV $n \in \mathcal{N}$
in time slot $t \in \mathcal{T}$ for a specific location ($A$, $B$, or $C$), respectively, as given by
\begin{subequations}
\begin{align}
    & x_{n, t} =  c_{n, t}^{A} - d_{n, t}^{A} \\
    & y_{n, t} =  c_{n, t}^{B} - d_{n, t}^{B} \\
    & z_{n, t} =  c_{n, t}^{C} - d_{n, t}^{C} \\
    &0 \leq c_{n, t}^{A} \leq P_{c}^{\max} \alpha_{n, t} \label{eq: A charge} \\
    &0 \leq d_{n, t}^{A} \leq P_{d}^{\max} \alpha_{n, t} \\
    &0 \leq c_{n, t}^{B} \leq P_{c}^{\max} \beta_{n, t} \\
    &0 \leq d_{n, t}^{B} \leq P_{d}^{\max} \beta_{n, t} \\
    &0 \leq c_{n, t}^{C} \leq P_{c}^{\max} \gamma_{n, t} \\
    &0 \leq d_{n, t}^{C} \leq P_{d}^{\max} \gamma_{n, t} \label{eq: C discharge} 
\end{align}
\end{subequations}
where vectors $\bm{\alpha}$, $\bm{\beta}$ and $\bm{\gamma}$ $\in \{0, 1\}$ are comprised of binary elements and indicate whether or not the EV $n \in \mathcal{N}$ is at
locations $A$, $B$ or $C$ at a time $t \in \mathcal{T}$, respectively. Constraints (\ref{eq: A charge})-(\ref{eq: C discharge}) describe the charging and discharging limits at locations $A$, $B$ or $C$, respectively, given by the available
bi-directional charger at those locations.\newline
Denoting the state of charge level of each EV $n \in \mathcal{N}$ in time slot $t \in \mathcal{T}$ by $E_n(t)$, the EV’s dynamics can be updated by
\begin{align}
    &\begin{split}
        &E_n(t) = E_n(t - 1) + \Big(\eta_n^c(c_{n, t}^{A}+c_{n, t}^{B}+c_{n, t}^{B})-\\&\frac{d_{n, t}^{A}+d_{n, t}^{B}+d_{n, t}^{C}}{\eta_n^d} - P^{drive}_n(1-\alpha_{n, t}-\beta_{n, t}-\gamma_{n, t})\Big) \Delta t
    \end{split}
\end{align}
where $E_n(0) = E_n^{\text{init}}$ and $P^{drive}_n$ represent the power consumption rate of EV while driving at each time step of $\Delta t$. The SOC of each EV $n \in \mathcal{N}$ in the time slot $t \in \mathcal{T}$ is bounded by the lower and upper SOC limit, denoted by $E^{\min}$ and $E^{\max}$, respectively. We have
\begin{align}\label{eq: SOC limit}
    E^{\min}_n \leq E_n(t) \leq E^{\max}_n
\end{align}
Constraints (\ref{eq: ch/dis variable})-(\ref{eq: SOC limit}) guarantee that the battery of each EV is within their physical range when providing delivery or doing arbitrage.
Additionally, to safeguard against the complete depletion of the electric vehicles' batteries at the close of each day, our model includes a constraint guaranteeing a sufficient state of charge, thereby ensuring the uninterrupted execution of the delivery tasks, denoted as \( E_n(T) = E_n^{\text{final}} \).

\subsubsection{Spatial-temporal Model}
We identify two possible states of the considered EVs: i)  parking EVs at the charging stations for park-and-charge, and ii) driving on the road and going to the delivery points. As mentioned earlier, the present state of an EV $n \in \mathcal{N}$ in the time slot $t \in \mathcal{T}$ at locations $A, B$, and $C$ are described by binary variables $\bm{\alpha}$, $\bm{\beta}$ and $\bm{\gamma}$, respectively. Moreover, we use $T_{A, B}$, $T_{A, C}$ and $T_{B, C}$ to define the travel/driving time from regions $A \leftrightarrow B$, $A \leftrightarrow C$, and $B \leftrightarrow C$, respectively. Therefore, the spatial-temporal constraints of the EVs can be expressed as follows:
\begin{subequations}\label{eq: spatial-temporal eq}
\begin{align}
    &\alpha_{n, t}+\beta_{n, t}+\gamma_{n, t} \leq 1 \label{eq: all binary}\\
    &\alpha_{n, t} + \beta_{n, t+\tau} \leq 1 \qquad \text{if} \quad \alpha_{n, t} =1  \qquad \tau \in [1,T_{A,B}] \label{eq: time travel A}\\
    &\beta_{n, t} + \alpha_{n, t+\tau} \leq 1 \qquad \text{if} \quad \beta_{n, t} =1  \qquad \tau \in [1,T_{A,B}]\\
    % &\alpha_{n, t-\tau} + \beta_{n, t} \leq 1  \qquad \tau \in [1,T_{A,B}]\\
    % &\beta_{n, t-\tau} + \alpha_{n, t} \leq 1  \qquad \tau \in [1,T_{A,B}]\\
    &\alpha_{n, t} + \gamma_{n, t+\tau} \leq 1 \qquad \text{if} \quad \alpha_{n, t} =1 \qquad \tau \in [1,T_{A,C}]\\
    &\gamma_{n, t} + \alpha_{n, t+\tau} \leq 1 \qquad \text{if} \quad \gamma_{n, t} =1 \qquad \tau \in [1,T_{A,C}]\\
    &\beta_{n, t} + \gamma_{n, t+\tau} \leq 1  \qquad \text{if} \quad \beta_{n, t} =1 \qquad \tau \in [1,T_{B,C}]\\
    &\gamma_{n, t} + \beta_{n, t+\tau} \leq 1 \qquad \text{if} \quad \gamma_{n, t} =1 \qquad \tau \in [1,T_{B,C}] \label{eq: time travel C}\\
    % &\alpha_{n, t-\tau} + \gamma_{n, t} \leq 1  \qquad \tau \in [1,T_{A,C}]\\
    % &\gamma_{n, t-\tau} + \alpha_{n, t} \leq 1  \qquad \tau \in [1,T_{A,C}]\\
    % &\beta_{n, t-\tau} + \gamma_{n, t} \leq 1  \qquad \tau \in [1,T_{B,C}]\\
    % &\gamma_{n, t-\tau} + \beta_{n, t} \leq 1  \qquad \tau \in [1,T_{B,C}] \label{eq: time travel C} \\
    & a_n \geq \alpha_{n, t}, \qquad a_n \leq \sum_{t \in \mathcal{T}} \alpha_{n, t}\\
    & b_n \geq \beta_{n, t}, \qquad b_n \leq \sum_{t \in \mathcal{T}} \beta_{n, t}\\
    & c_n \geq \gamma_{n, t}, \qquad c_n \leq \sum_{t \in \mathcal{T}} \gamma_{n, t}\\
    & \sum_{n \in \mathcal{N}} a_n \geq N_A \label{eq: delivery A}\\
    & \sum_{n \in \mathcal{N}} b_n \geq N_B\\
    & \sum_{n \in \mathcal{N}} c_n \geq N_C \label{eq: delivery C}\\
    % &  \sum_{t \in \mathcal{T}} \sum_{n \in \mathcal{N}} \alpha_{n, t} = N_A \label{eq: delivery A}\\
    % &  \sum_{t \in \mathcal{T}} \sum_{n \in \mathcal{N}} \beta_{n, t} = N_B\\
    % &  \sum_{t \in \mathcal{T}} \sum_{n \in \mathcal{N}} \gamma_{n, t} = N_C \label{eq: delivery C} \\
    & \alpha_{n, 0} = 1 \qquad  n \in \mathcal{J}_{A} \label{eq: start loc A}\\
    &\beta_{n, 0} = 0 \qquad  n \in \mathcal{J}_{B} \\
    &\gamma_{n, 0} = 0 \qquad  n \in \mathcal{J}_{C} \label{eq: start loc C}
\end{align} 
\end{subequations}
It is not possible for each vehicle to be present at points $A, B$, and $C$ at the same time. Therefore, constraint (\ref{eq: all binary}) ensures that the vehicle $n \in \mathcal{N}$ is not at points $A, B$, and $C$ in the time slot $t \in \mathcal{T}$ simultaneously. While the EV is driving between two locations, all the binary variables should remain zero until the EV reaches the destination. constraints (\ref{eq: time travel A}) - (\ref{eq: time travel C}) aims to ensure that binary variables of the associated locations remain zero until the time steps after the vehicle $n \in \mathcal{N}$ has departed the first location. For the delivery option, it assumed that at least a specific number of vehicles should visit each location (constraints (\ref{eq: delivery A})-(\ref{eq: delivery C})).
A subset of vehicles is assumed to start at locations $A, B$, and $C$ in the first timestep, respectively, which necessitate constraints (\ref{eq: start loc A})-(\ref{eq: start loc C}).

\subsubsection{Stochastic Problem Formulation}
We have formulated the optimization problem under the assumption that the arrival time, departure time, and charging demand of electric vehicles at each location are unknown. We aim to find the optimal charging solution for each EV $n$ to minimize the expected total EV charging and discharging cost in an operational horizon $T$ (here, 24 hours with 15-minute time steps). In other words, this objective seeks to maximize the expected profits derived from spatial and temporal arbitrage, while taking into account the associated expenses of the EV charging and discharging processes. Then we can formulate the charging optimization problem for the aggregator owing fleet of EVs as follows.
\begin{align}
\begin{split}  
    &\underset{\substack{\bm{x}, \bm{y}, \bm{z}}}{\text{Minimize}} \quad  \mathbb{E}[\bm{x}, \bm{y}, \bm{z}, \bm{\xi}_t^A, \bm{\xi}_t^B, \bm{\xi}_t^C]  = \\ &\frac{1}{N_k}\sum_{k \in \mathcal{K}} \sum_{t \in \mathcal{T}} \sum_{n \in \mathcal{N}} (\xi^A_{k, t} x_{n, t} + \xi^B_{k, t} y_{n, t}+\xi^C_{k, t} z_{n, t})
\end{split}\\
&\textrm{subject to:} \nonumber\\
&\qquad \textrm{Constraints (\ref{eq: ch/dis variable})-(\ref{eq: spatial-temporal eq})} \nonumber
\end{align}
where the decision variables are within vectors $\bm{x}, \bm{y}, \bm{z}, \bm{\alpha}, \bm{\beta}$, and $\bm{\gamma}$. Here, uncertain parameters are the electricity price at each time $t \in \mathcal{T}$ at location $A$, $B$, and $C$, $\bm{\xi}_t^A$, $\bm{\xi}_t^B$ $\bm{\xi}_t^C$, respectively. These uncertain parameters are assumed to have a finite number of considered realizations $|\mathcal{K}|$ and equal probability masses, although this can be modified if one knows that certain scenarios are more likely
to be similar to the present situation, for example.

% \begin{align}
% 	\mathcal{L}(\bm{\theta}; \mathcal{T}) = w_f\mathcal{L}_f(\bm{\theta}; \mathcal{T}_f) + w_b\mathcal{L}_b(\bm{\theta}; \mathcal{T}_b)\label{eq: PINN-loss}, 
% \end{align}where
% \begin{align*}
% \begin{split}
%     &\mathcal{L}_f(\bm{\theta} ; \mathcal{T}_f) = \frac{1}{|\mathcal{T}_f|}\sum_{(\bm{x}, t) \in \mathcal{T}_f} |f(\bm{x}, t)|^2\\ 
%     &\mathcal{L}_b(\bm{\theta}; \mathcal{T}_b) =\frac{1}{|\mathcal{T}_b|}\sum_{(\bm{x}, t) \in \mathcal{T}_b} |\mathcal{B}( \hat{\bm{u}}, \bm{x})|^2
% \end{split}
% \end{align*}

% \section{Physical Model for Power System Dynamics} \label{sec: model}

\section{Case Study}
\begin{figure}[t!]
\centering
\includegraphics[clip,trim={0.3cm 1.35cm 2cm 0.7cm},width=\columnwidth]{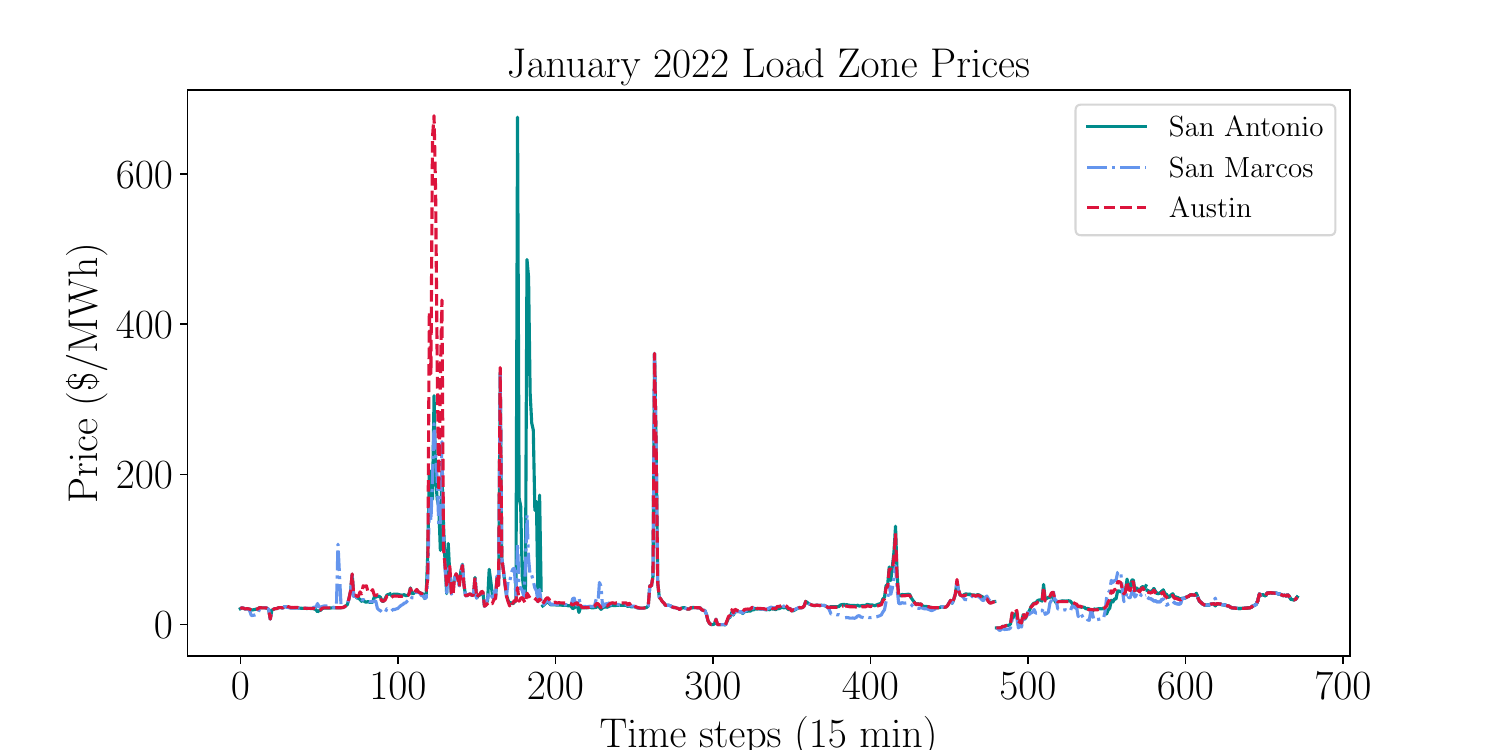}\\
\includegraphics[clip,trim={0cm 0cm 2cm 0cm},width=\columnwidth]{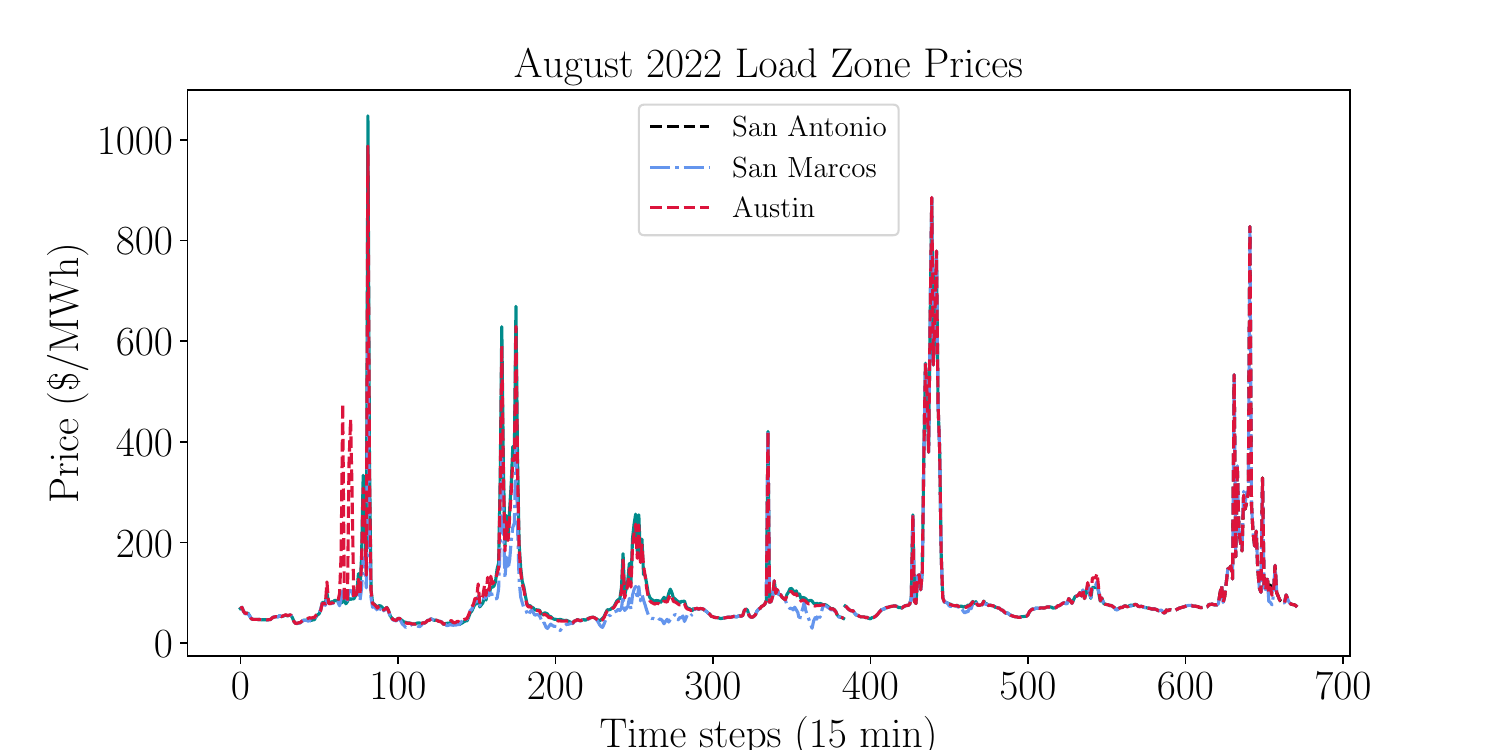}\\%
\caption{Historical real-time market settlement point prices for January and August 2022 in the three geographical locations. }
\label{fig:price-Jan}
\end{figure}

% We consider multiple case studies to test the proposed framework.
% \subsection{Simulation Setup and Training}
In our study, we utilized two weeks from January and August $2022$ for simulations. We sourced $15$-minute interval real-time zonal pricing data from the Electric Reliability Council of Texas (ERCOT) load zones to determine electricity prices for San Antonio, San Marcos, and Austin, as depicted in Fig. \ref{fig:price-Jan} for the associated periods. This visual representation delineates the fluctuating electricity prices across three prominent zones: San Antonio (solid line), San Marcos (dash-dotted line), and Austin (dashed line). Notably, sharp price spikes are evident at specific intervals, most pronounced in the time steps between $100$ and $300$. While the majority of the data showcases relatively stable price ranges, these intermittent spikes underscore the inherent volatility and complexity of real-time electricity pricing in the Texas grid. These intermittent peaks will be relevant when later analyzing the arbitrage strategy employed.

Given that a fleet of EVs doesn't represent a significant energy resource, we posited that the electric delivery trucks would transact energy based on zonal prices, and the vehicles are price-takers (e.g., they are unable to directly influence market electricity prices through their charging/discharging actions). These prices represent an average across all nodes in a particular zone. In contrast, nodal prices, which are generally more variable and are typically employed by larger, centralized resources, were not considered for this analysis. Note, however, that these prices could easily be incorporated into the framework in place of the zonal assumption.

\begin{figure}[t!]
\centering
\includegraphics[width=\columnwidth]{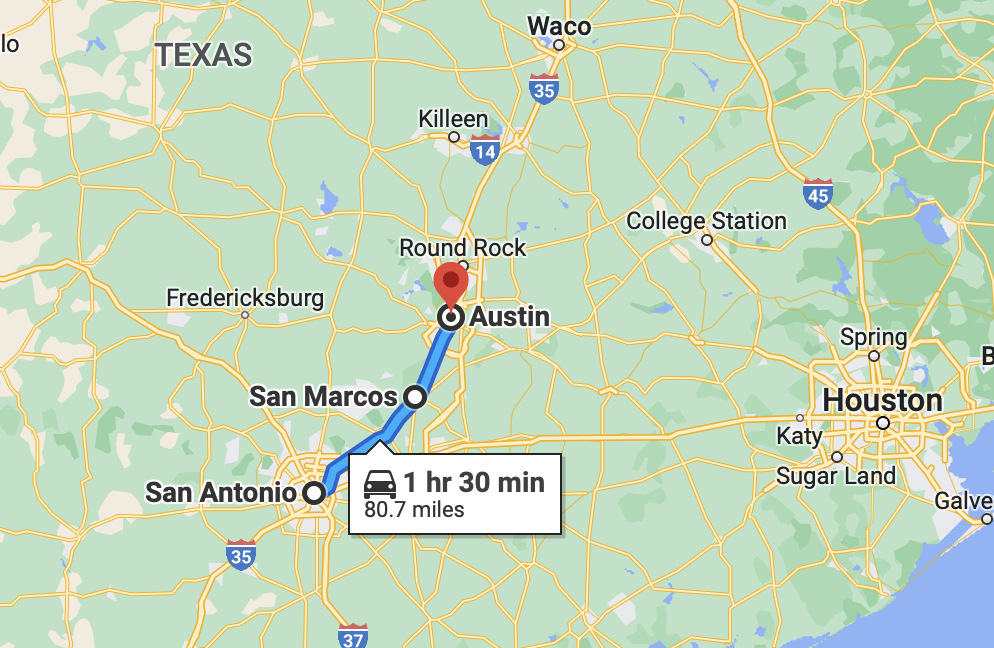}\\%
\caption{Relative location of the three considered cities within Texas. Typical traffic for each hour and day of the week was obtained from Google Maps.}
\label{fig:traffic}\vspace{-4mm}
\end{figure}

Due to the unavailability of public historical traffic data for the targeted regions, we relied on fifteen-minute ``typical traffic'' data for weekdays, sourced from Google Maps. This data is derived from historical traffic patterns to estimate travel times between the cities. Fig. \ref{fig:traffic} illustrates the relative locations of San Antonio, San Marcos, and Austin: roughly 31 miles separate Austin and San Marcos, and 50 miles separate San Marcos to San Antonio. This provides a reasonable proxy but may not fully capture real-world congestion and travel time variability.

For the present simulations, we modeled electric delivery trucks equipped with a battery capacity ranging from $630$ to $770$ kWh. This capacity was randomly assigned to each truck. Each battery also has a random roundtrip efficiency between $90-100\%$.  These trucks have access to DC Fast Charging stations with a capacity of $150$ kW. We operated under the assumption that these charging stations could facilitate bidirectional charging at the same rate, although the framework can incorporate different charge/discharge limits. Throughout the week, if an EV opted for spatial arbitrage, it would begin the subsequent day from its starting location.
In our base scenario, we considered a fleet of 10 delivery trucks, with each truck making a minimum of six deliveries daily. The initial state of charge for these EVs was randomized between $420$ to $490$ kWh. Furthermore, the consumption rate while driving was estimated to be between $63$ to $77$ kWh per hour if an average of 60 miles/hour is assumed.
Lastly, our framework was developed in Python and employed the Gurobi solver \cite{Gurobi} for optimization. All simulations were executed locally on a personal Apple laptop equipped with the M1 chip.

% Central case study results

%\begin{figure}[h!]
%\centering
% \includegraphics[clip,trim={0cm 0cm 0cm 1cm},width=\columnwidth]{figs/IMG_1786.jpg}%
%\includegraphics[clip,trim={.95cm 0cm 0cm 1cm},width=\columnwidth]{figs/x_aug.pdf}\\\vspace{-5mm}%
%\includegraphics[clip,trim={.95cm 0cm 0cm 1cm},width=\columnwidth]{figs/y_aug.pdf}\\\vspace{-5mm}%
%\includegraphics[clip,trim={.95cm 0cm 0cm 1cm},width=\columnwidth]{figs/z_aug.pdf}%
%\caption{The cumulative charging (and discharging) of the fleet of delivery vehicles for each of the three locations. The electricity price at that location is superimposed onto each subplot.}
%\label{fig:fleet_cd}
%\end{figure}

\subsection{Considered Scenarios}

To understand the benefits of performing spatial arbitrage, we considered three different scenarios:

\begin{itemize}
    \item \emph{Spatial Arbitrage:} The delivery fleet must adhere to a given delivery schedule, can both charge and discharge, and the electricity price changes throughout time and space.
    \item \emph{Counterfactual:} The delivery fleet must adhere to a given delivery schedule, can both charge and discharge, but the delivery scheduling is does not take into account prices which vary by location; only by time.
    \item \emph{Stationary Arbitrage:} The delivery fleet is assumed to be stationary, but can charge and discharge at the warehouse location (San Antonio).
\end{itemize}

The purpose of the \emph{Counterfactual} case is to remove the incentive of the EVs to travel or take unnecessary trips due to prices that vary across location, yielding a result that adheres to the delivery constraints while still having the option to perform temporal arbitrage. The profits reported for this case are then calculated using the actual time varying prices.

Table \ref{tab: maincase} elucidates the economic and operational outcomes derived from a week-long simulation optimizing the use of a fleet of ten electric delivery trucks engaged in arbitrage across three scenarios. The scenarios are crafted to assess the impact of strategic movement and charging decisions on the cost-effectiveness, distance traveled, and energy consumption within the fleet operations. In January, the \emph{Spatial Arbitrage} scenario shows the highest profitability, with the fleet earning $\$3,996.63$. This represents a $42.7\%$ increase in profits compared to the \emph{Counterfactual} scenario, where the trucks earn $\$2,786.90$ despite adhering to the same number of deliveries. This comparison demonstrates the added value of spatial arbitrage when distinct electricity pricing is available across locations. The \emph{Stationary Arbitrage} scenario, where trucks do not travel, still manages a profit of $\$2,961.74$, which is $5.9\%$ higher than the \emph{Counterfactual} scenario but $25.9\%$ less than the \emph{Spatial Arbitrage} scenario, emphasizing the potential gains from engaging in spatial arbitrage even when the trucks are also tasked with deliveries.

Moving to August, the profitability for the \emph{Spatial Arbitrage} scenario increases dramatically to $\$7,851.94$, which is a $96.8\%$ surge from January's profits. This dramatic rise suggests a higher efficiency in operational strategy or potentially greater disparities in electricity prices during this month, leading to more lucrative arbitrage opportunities. In comparison, the \emph{Counterfactual} scenario's profit in August is $\$5,316.51$, a significant $90.8\%$ increase from January's earnings, indicating that even without spatial arbitrage, the seasonal variations in electricity prices or delivery schedules can significantly affect profitability. When we compare the two scenarios, the \emph{Spatial Arbitrage} profits in August outperform the \emph{Counterfactual} by $47.7\%$, which is a slightly higher margin than seen in January ($42.7\%$). This increase suggests that spatial arbitrage becomes even more advantageous when the differential in electricity pricing between locations widens.

In the \emph{Stationary Arbitrage} scenario, where trucks remain stationary and do not partake in delivery tasks, there is a purely financial focus on energy arbitrage. This case shows a consistent approach, with zero distance traveled and a slight decrease in profits from January to August by $5.8\%$, with a profit of $\$1,009.91$ in August. This decrease could be due to a lower potential for temporal arbitrage or an overall decrease in energy prices during this month, reflecting the varying dynamics of the energy market across different times of the year. In summary, spatial arbitrage consistently shows a higher profit margin over the stationary approach, highlighting the potential for significant earnings through strategic movements and energy trading by EV fleets. Therefore, the seasonal variations in profit margins point to the substantial impact that market factors can have on arbitrage strategies.

\begin{table*}[ht!]
    \centering 
    \vspace{2mm}
\caption{The cost, distance driven, and total battery throughput over the week-long simulation. We compare the spatial arbitrage case with the delivery counterfactual, the case where bi-directional charging is performed only at the warehouse.}
\label{tab: maincase}
\addtolength{\tabcolsep}{-0.15em}
\begin{tabular}{l c c c c c c}
\toprule
\multirow{2}{*}{Scenario}  & \multicolumn{3}{c}{\textbf{January}} & \multicolumn{3}{c}{\textbf{August}} \\  
\cmidrule(lr{0.5em}){2-4} \cmidrule(lr{0.5em}){5-7} 
    &  Cost ($\$$) & Distance (mi) & \makecell{Energy \\Consumption (kWh)} &  Cost ($\$$) &  Distance (mi) &  \makecell{Energy \\Consumption (kWh)}\\
\midrule
Spatial Arbitrage & $-3996.63$ & $10813.40$  & $87262.10$ & $-7851.94$ & $9000.30$  & $81349.74$\\
Counterfactual & $-2786.90$ & $6570.90$  & $76798.62$  & $-5316.51$ & $6726.09$  & $78098.00$\\
Stationary Arbitrage & $-2961.74$  & $0.00$  & $75611.50$ & $-1009.91$  & $0.00$  & $69151.15$ \\
\bottomrule
\end{tabular}
\vspace{-3mm}
\end{table*}

\begin{figure}
    \centering
    \includegraphics[width=\columnwidth]{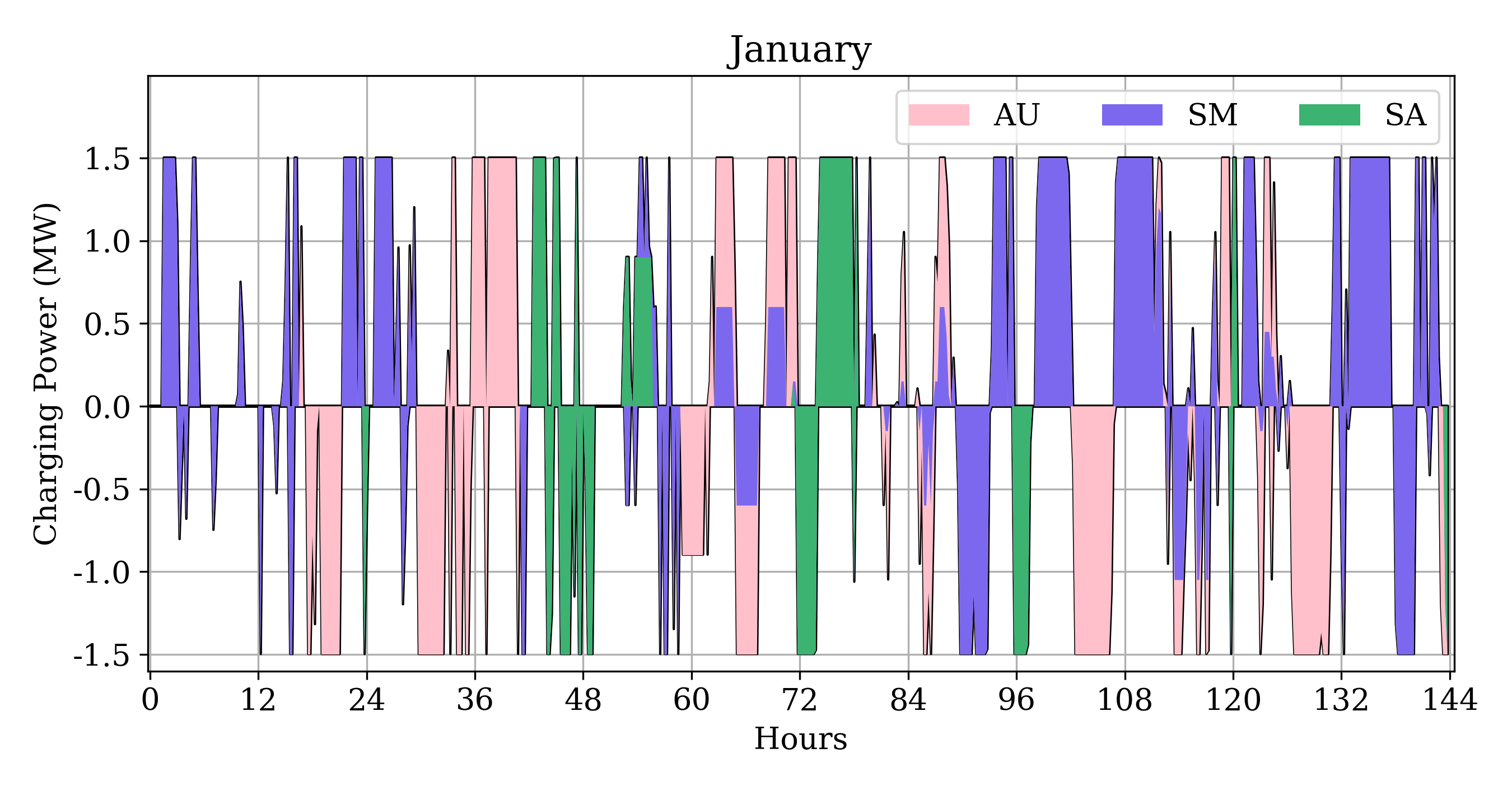}
    \includegraphics[width=\columnwidth]{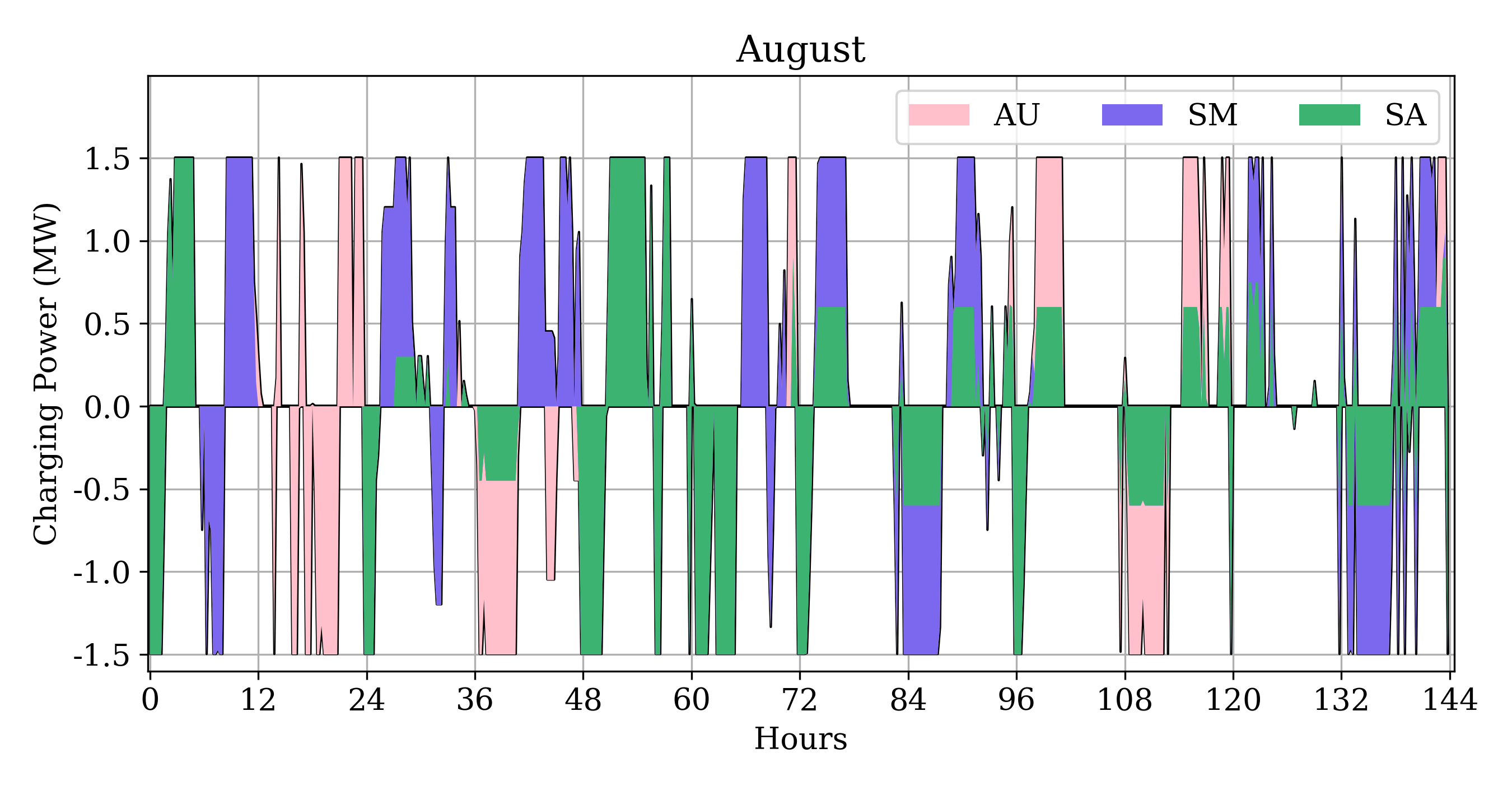}
    \caption{Cumulative fleet charging and discharging broken down by the three locations (AU = Austin, SM = San Marcos, SA = San Antonio), for a week in January and a week in August.}
    \label{fig:net_charging}
\end{figure}

Figure \ref{fig:net_charging} presents the cumulative net charging activities of a fleet of electric EVs engaged in spatial arbitrage, delineated by location, for one week in January and one week in August. The vertical axis represents the cumulative power (in MW), while the horizontal axis denotes the hours over the week. Note that positive cumulative power relates to battery charging activities, while negative cumulative power indicates discharging back to the grid. The color coding distinguishes the locations: pink for Austin (AU), purple for San Marcos (SM), and green for San Antonio (SA). In both months, the fleet's charging (positive values) and discharging (negative values) behaviors demonstrate a clear response to the price signals indicated in Fig. \ref{fig:price-Jan}, which portrays the historical real-time market settlement point prices for the same weeks in January and August. The strategic charging and discharging, i.e., energy arbitrage, capitalize on the fluctuating electricity prices to generate profit, charging the batteries when prices are low and discharging when prices are high.

August displays an intensified pattern of arbitrage activities, reflecting the higher price volatility illustrated in Fig. \ref{fig:price-Jan}. The escalated pricing dynamics in August, particularly noticeable with Austin's steep and frequent price peaks, correlate with the larger profits achieved in the \emph{Spatial Arbitrage} scenario for this month, as outlined in the earlier provided Table \ref{tab: maincase}. 
% The negative cost values indicate profits, and the larger absolute value for August (-$7851.94) compared to January (-$3996.63) highlights greater profitability during this summer month.
The consistent pattern across both months confirms that the fleet of EVs is utilized not only for energy storage but also as a mobile energy resource that shifts energy geographically according to arbitrage opportunities. This is in line with the optimization problem's objective to maximize profits through strategic spatial and temporal energy trading while adhering to physical constraints such as a vehicle's presence in a single location at any given time. Such integration could contribute to stabilizing grid prices and providing a buffer against extreme price fluctuations.

\begin{figure}[t!]
\centering
\includegraphics[clip,trim={0.3cm 1.35cm 0cm 0.7cm},width=\columnwidth]{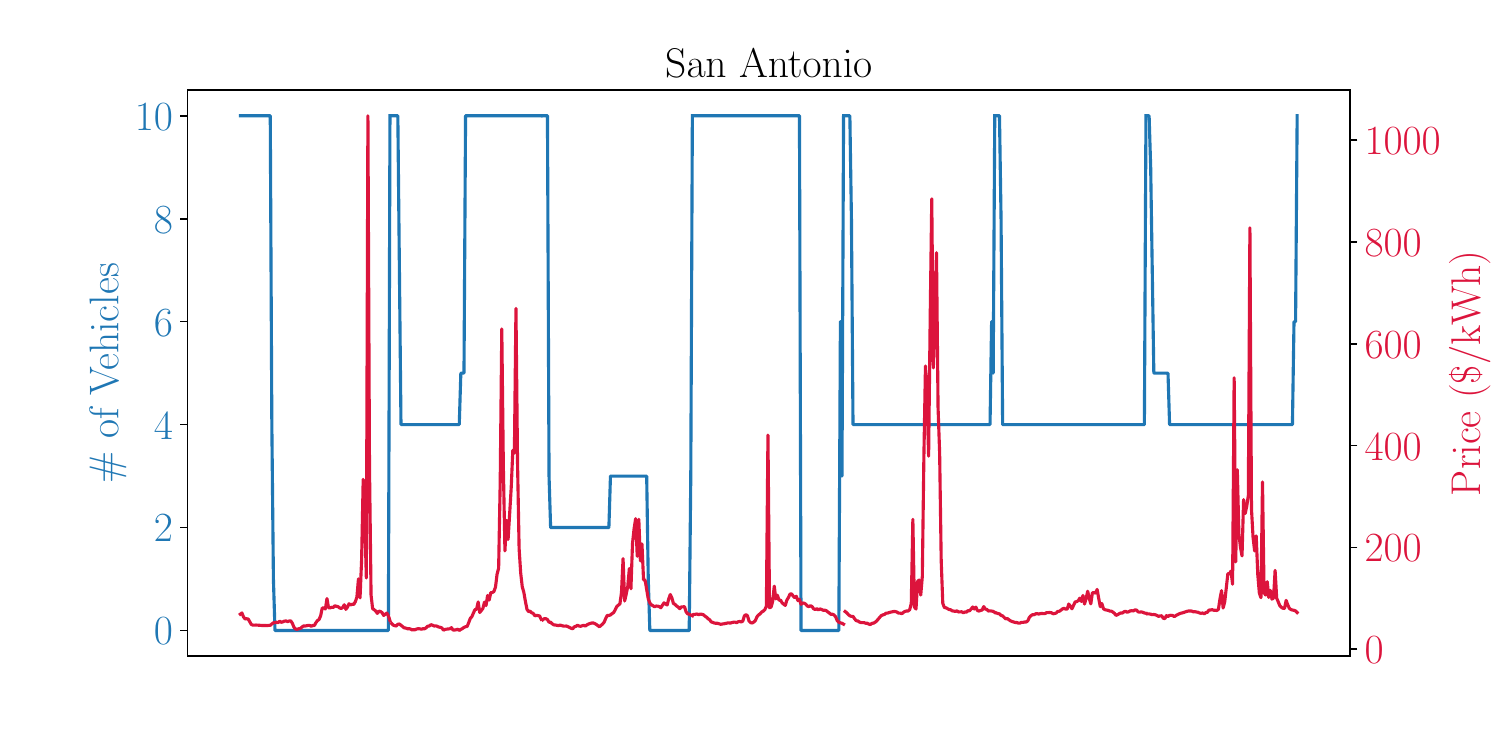}\\
\includegraphics[clip,trim={0.3cm 1.35cm 0cm 0.7cm},width=\columnwidth]{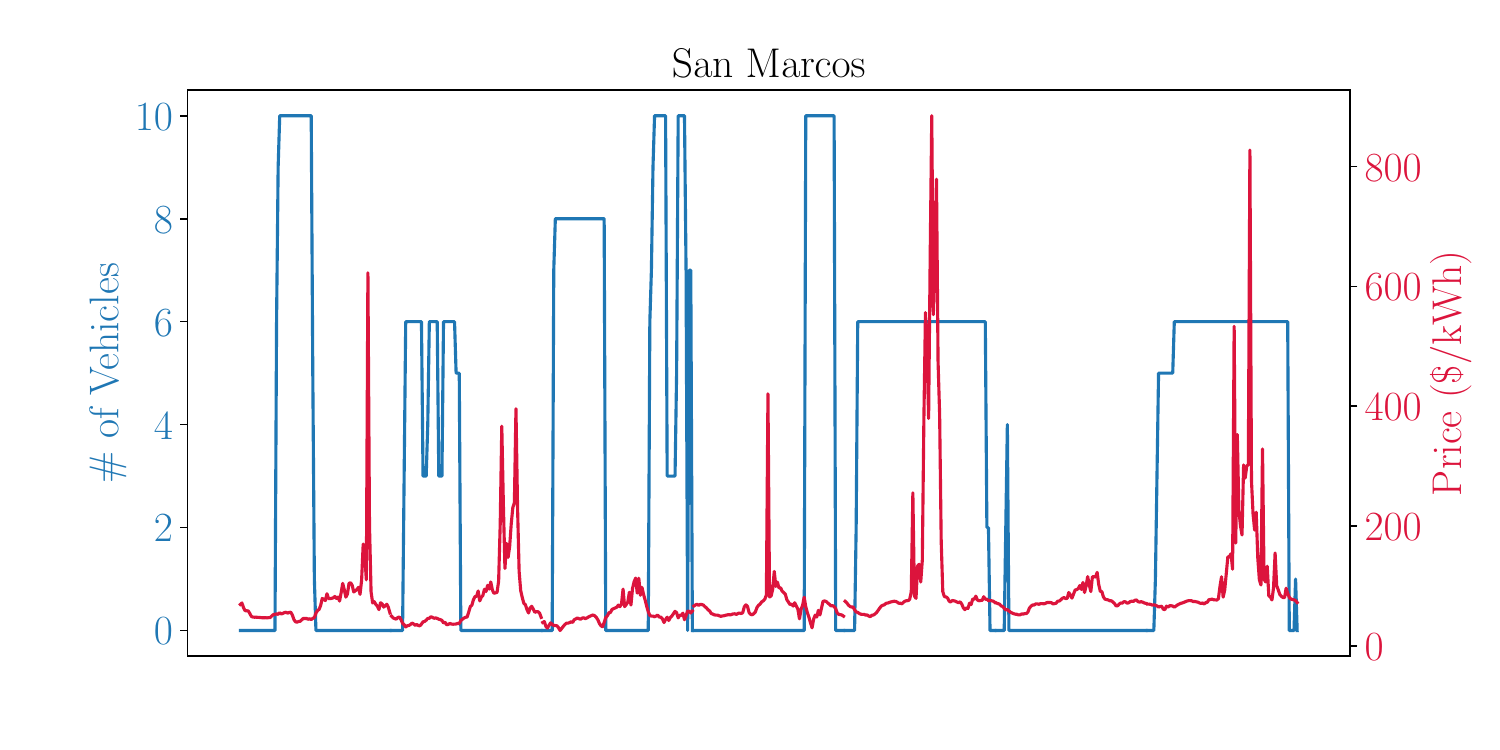}\\
\includegraphics[clip,trim={0cm 0cm 0cm 0.7cm},width=\columnwidth]{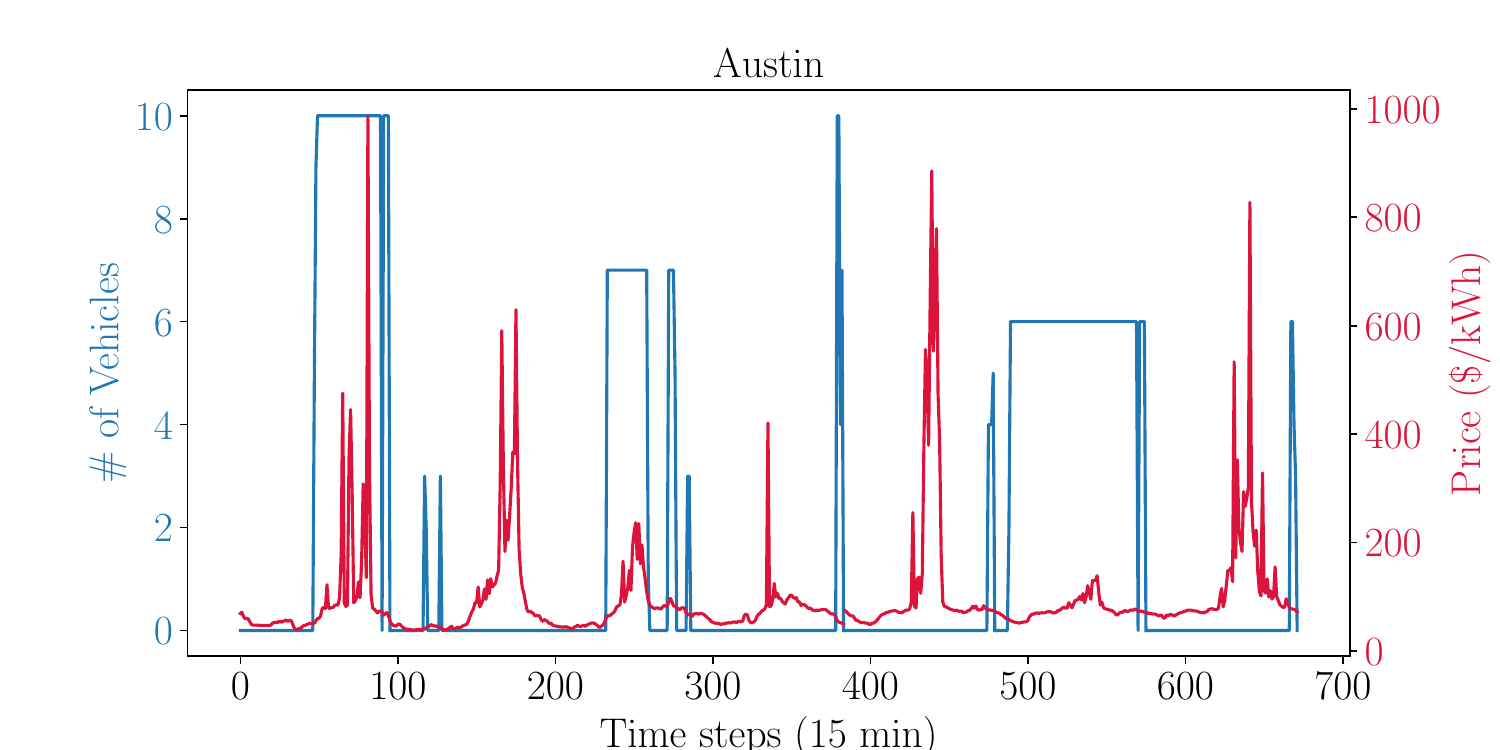}\\%
\caption{Distribution and pricing dynamics of the EV fleet performing spatial arbitrage in August}
\label{fig: Num vehicles}
\end{figure}

Figure \ref{fig: Num vehicles} illustrates the distribution of a fleet of electric vehicles across San Antonio, San Marcos, and Austin, juxtaposed with the corresponding real-time electricity market prices for each location. It captures the fleet's strategic positioning throughout a week in August, reflecting the balancing act between fulfilling delivery requirements and capitalizing on arbitrage opportunities by leveraging temporal and spatial price differentials in the energy market.
In San Antonio, the fluctuation in the number of vehicles corresponds with shifts in energy pricing, suggesting that vehicle distribution is, at least in part, a reaction to economic incentives. During peak price periods, the presence of fewer vehicles may reflect a strategic dispatch to other locations, likely in pursuit of lower charging costs or higher revenue from discharging. Conversely, a rise in the number of vehicles during off-peak price periods suggests a collective strategy to capitalize on cheaper charging opportunities.

The San Marcos graph indicates a similar, though not as pronounced, pattern of vehicle distribution in response to price changes. This is likely a result of the geographic and economic positioning of San Marcos between San Antonio and Austin, making it a strategic midpoint for charging and discharging activities, as well as fulfilling the minimum delivery requirement. In Austin, the pronounced spikes in energy prices align with substantial drops in the number of vehicles, which supports the notion that vehicles are actively being repositioned in anticipation of or in response to these price surges. This behavior underscores the sophisticated optimization model at play, where the delivery schedule is maintained (as evidenced by at least six vehicles consistently being dispatched for deliveries) while still maximizing arbitrage profits. The results in Table \ref{tab: maincase} highlighted the profitability of such operations, with August showing a more substantial profit margin than January. This figure underpins those results, illustrating how the fleet's dynamic redistribution aligns with economic incentives to optimize profits while adhering to delivery obligations. This fleet's ability not only reinforces the financial viability of such operations but also highlights the potential of EV fleets as flexible energy resources that can support grid stability.

\begin{figure}[h!]
\centering
\includegraphics[clip,trim={0cm 0cm 0cm 0cm},width=\columnwidth]{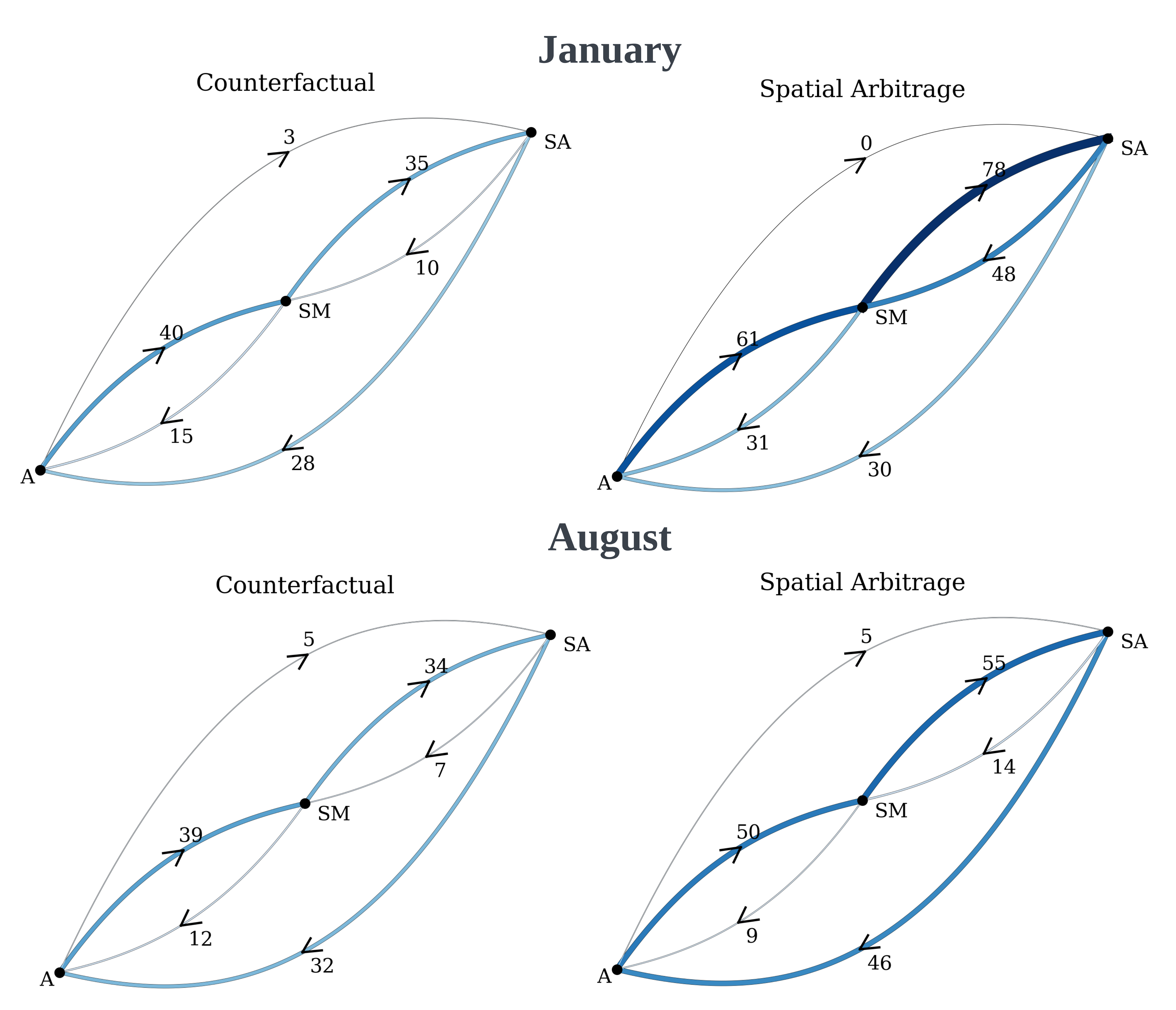}%
\caption{A visualization of the number of trips carried out in the spatial arbitrage scenario versus the counterfactual delivery schedule, for both January and August. Note that the arrow widths are proportional to the number of trips.}
\label{fig:map}
\end{figure}

Figure \ref{fig:map} provides a comparative analysis of the number of trips made by a fleet of EVs under two scenarios: a counterfactual delivery schedule and a spatial arbitrage strategy during one-week periods in January and August. The diagram illustrates the directional flow and volume of trips between three locations—Austin (A), San Marcos (SM), and San Antonio (SA)—with the thickness of the arrows indicating the relative number of trips. 
Through a detailed examination, it is apparent that in January, the spatial arbitrage model incurs a higher number of trips compared to the counterfactual, yet the profits do not proportionally reflect this increase due to weaker price correlations. For instance, while the trip count increases by $50\%$ from the counterfactual to the arbitrage scenario, profits may not exhibit a commensurate rise. Conversely, August shows a more strategic engagement with fewer trips but significantly greater profitability, suggesting that arbitrage efficiency is highly dependent on market conditions and price volatility.

Energy consumption from additional trips is overshadowed by the energy involved in arbitrage operations. In other words, the adoption of spatial arbitrage is accompanied by an increase in energy usage, predominantly driven by heightened battery cycling for charging and discharging. Nevertheless, this increase in energy consumption, which could be hypothetically pegged at $30\%$, is compensated by a more pronounced surge in profits, perhaps exceeding $50\%$, highlighting the method's effectiveness despite its energy demands. Despite the lower absolute profit levels in January, spatial arbitrage still represents a considerable improvement over the counterfactual approach. If the profits in the counterfactual scenario are considered a baseline, spatial arbitrage could represent an increase to $150\%$ of this baseline in January and up to $200\%$ in August.

The analysis also reveals that the additional energy consumed for the trips mandated by spatial arbitrage is substantially less than the overall energy consumption recorded. This underscores that while the baseline energy expenditure for the minimum required trips is comparatively modest, the actual consumption under arbitrage is higher, reflecting the trade-offs for increased profit margins. The variance in the number of trips between January and August can be attributed to the better alignment of electricity prices between locations in the latter month, which enables a more effective arbitrage strategy with reduced trip requirements to maintain or enhance profitability. Additionally, it is crucial to recognize that the mileage accrued under spatial arbitrage is likely to be greater than that under the counterfactual. This must be weighed against the consequent implications for vehicle maintenance, wear-and-tear, and the longevity of the battery.

Despite January's intensive operational demands, spatial arbitrage displays a definitive advantage in terms of profitability and operational efficiency, especially pronounced in August. This analysis furnishes valuable insights into the feasibility and potential of spatial arbitrage as a strategic operational approach for EV fleets. By exploiting the temporal and spatial fluctuations in energy prices, fleet operators have the potential to significantly elevate profit margins. However, this must be prudently balanced against energy consumption concerns, battery health, and the overarching aim of sustainable operations.

% The spatial arbitrage scenario indicates strategic repositioning of the fleet to exploit electricity price differentials across locations, whereas the counterfactual represents standard delivery routing without considering energy market opportunities.

% Energy consumption from additional trips is overshadowed by the energy involved in arbitrage operations. However, the increase in profits outstrips the rise in energy consumption in percentage terms for both months.
% The observed profits in the spatial arbitrage scenario for both months are significantly higher than those achieved in the counterfactual delivery schedule, emphasizing the effectiveness of the arbitrage strategy.

\begin{figure}[t!]
\centering
\includegraphics[clip,trim={0.3cm 0cm 0cm 0.7cm},width=\columnwidth]{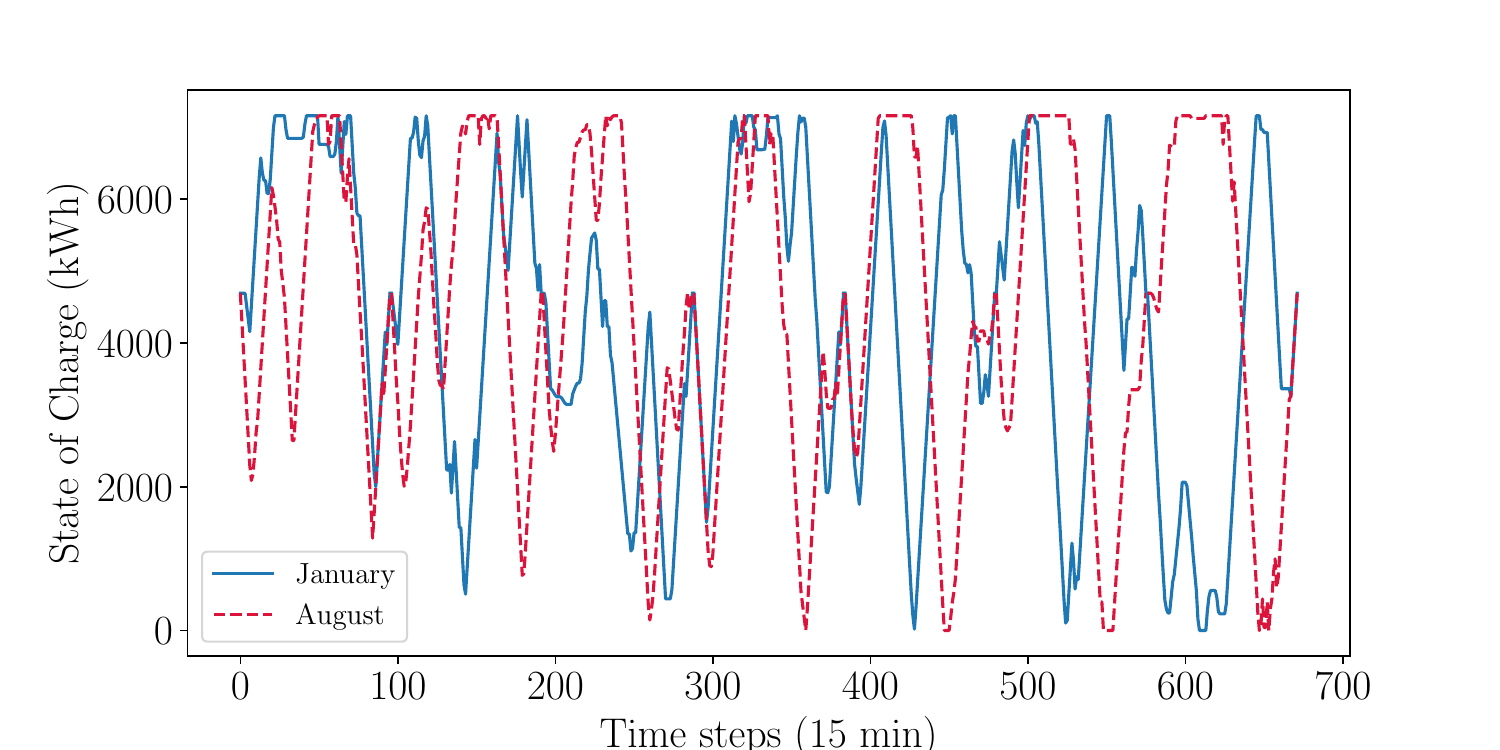}\\
\caption{Cumulative state of charge of the EVs across spatial arbitrage scenario.}
\label{fig: SOC}
\end{figure}

The cumulative state of charge (SOC) of EVs operating under a spatial arbitrage scenario for January and August is displayed in Fig. \ref{fig: SOC}. It illustrates the charging and discharging cycles of the EVs corresponding to the energy arbitrage opportunities that arise due to price fluctuations.

In August, the SOC pattern shows more pronounced fluctuations compared to January, suggesting a more active engagement in arbitrage as the vehicles frequently adjust their SOC in response to changing energy prices. This indicates a strategy that aggressively capitalizes on higher price differentials during this month. The graph also shows that while the overall SOC trends follow a similar pattern across the two months, the amplitude and frequency of the fluctuations differ. This difference can be tied to the seasonality in electricity price volatility and demand for energy, which affects how often and how significantly EVs need to be charged or discharged within the spatial arbitrage model.

The requirement that SOC levels return to a baseline ($E_n^{\text{final}}$) at the end of each day (every 96-time steps) introduces a cyclic pattern in the graph and ensures that vehicles maintain enough charge to perform necessary deliveries without the risk of battery depletion. This constraint not only safeguards the operational reliability of the EVs but also reflects a balanced approach to energy management, aligning profit-making from arbitrage with the logistical demands of delivery schedules.

\subsection{Effect of uncertain prices}

\begin{figure}
    \centering
    \includegraphics[width=\columnwidth]{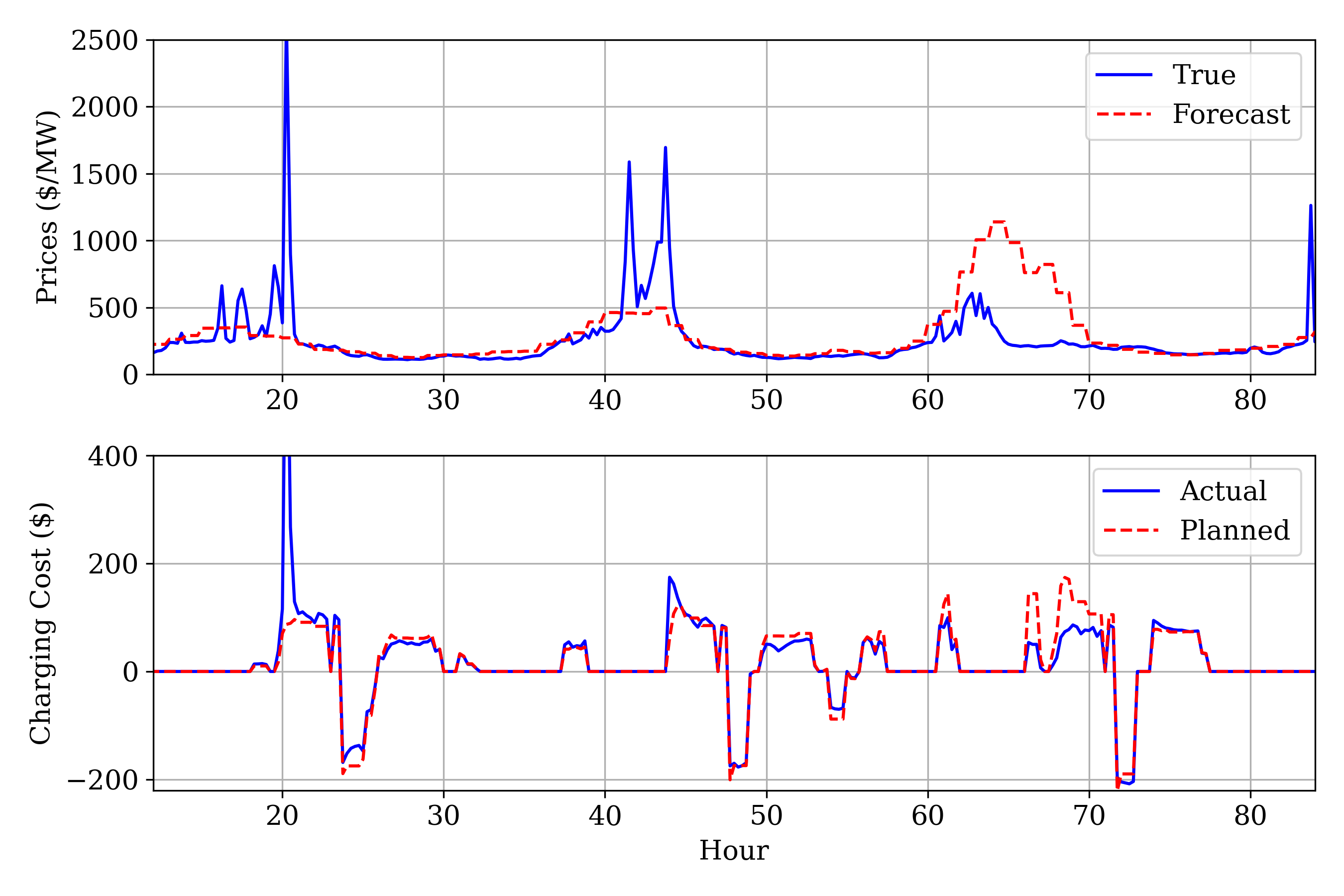}
    \caption{The performance of the planning algorithm in the August simulation when day-ahead prices are used as forecasts instead of assuming perfect knowledge of the real-time price.}
    \label{fig:uncertainty}
\end{figure}

In practice, it is not realistic to assume a perfect knowledge of the real-time prices. In this section, we analyze the effect of using the day-ahead price to plan charging/discharging while taking the real-time price. The August case study was repeated using the day-ahead prices for each region. Figure \ref{fig:uncertainty} shows both the day-ahead and real-time prices averaged over the three load zones, as well as the intended vs. actual fleet charging cost over the same time period. Day-ahead prices are less volatile than real-time but capture the general shape reasonably well. However, noticeably the day-ahead prices totally miss the extreme peaks in real-time price. This can adversely affect the fleet charging planning; on the first day shown the true fleet charging cost explodes for one time step due to a surprise peak in price. Despite this, in this simulation, there were also times when the charging costs were lower than planned. 

On average, the error in forecasting added around $\$25$/day in charging costs. For more severe unforeseen peaks, these costs could be worse. However, it is likely that the addition of a rules-based controller could produce a more conservative strategy -- for example one which stops all charging if the real-time price is too much higher than the forecast price. 

\section{Conclusion}\label{sec: conclude}
We introduced a streamlined computational framework tailored for a fleet of electric trucks, facilitating both temporal and spatial arbitrage. Demonstrated findings reveal that fleet-wide EV optimization, when managed by an aggregator distinct from delivery duties, enriches the financial portfolio's diversity. 
% This paper introduced a novel optimization framework to maximize the profitability of a fleet of electric delivery trucks through spatial and temporal energy arbitrage. 
The proposed model integrates charging, discharging, routing, and scheduling decisions into a single-stage stochastic program. Through simulations based on real-world pricing data, the results demonstrate the financial viability and operational feasibility of leveraging electric vehicle fleets as mobile energy storage to capitalize on locational marginal price differentials.

The case study analysis highlights several key findings. First, spatial arbitrage consistently outperforms both the counterfactual and stationary arbitrage benchmarks, with profits increasing by as much as $47.7\%$ compared to non-spatial optimization strategies. This substantiates the value derived from strategic positioning and energy trading across locations. Second, although additional trips are required to enable spatial arbitrage, the energy expenditures from these trips are marginal compared to consumption from charging and discharging. Therefore, the energy overhead is outweighed by profit gains. Third, seasonal pricing dynamics significantly impact arbitrage outcomes, with wider spreads and higher volatility enabling greater returns. This underscores the importance of real-time market awareness and forecasting.

While the results showcase promising potential, for subsequent research, it would be invaluable to factor in battery degradation costs, charging infrastructure limitations, and forecasting uncertainties as they play a pivotal role in financial deliberations. Constraints on battery cycling could possibly be included in the optimization to introduce a cost-benefit trade-off analysis. More complex constraints on delivery windows, locations, and large fleet management, in general, would also add to the fidelity of the simulation framework. 
Overall, the model provides a valuable framework to optimize electric vehicle operations at the nexus of transportation and energy system needs. This integration of mobility and electricity markets via vehicle-to-grid capabilities offers environmental and economic value for numerous stakeholders.
\bibliography{references.bib}{}
\bibliographystyle{IEEEtran}

%That's all
\end{document}